\documentclass[11pt,leqno,a4paper]{article}

\usepackage{amsmath, amssymb, amsthm, amsfonts}
\usepackage{mathrsfs}

\theoremstyle{plain}
\newtheorem{thm}{Theorem}

\newtheorem{lem}{Lemma}

\theoremstyle{remark}
\newtheorem{rem}{Remark}

\def\K{\mathop{\mbox{\bf\Large K}}}

\numberwithin{equation}{section}

\begin{document}

\title{Multiple-correction and continued fraction approximation(II)}
\author{ Xiaodong Cao and Xu You}
\date{Oct 10, 2014}

\maketitle



\footnote[0]{This work is supported by the National Natural Science Foundation of China (Grant No.11171344 and 61403034)
and the Natural Science Foundation of Beijing (Grant No.1112010).}

\begin{abstract}
The main aim of this paper is to further develop the multiple-correction method that formulated in our previous works~\cite{CXY, Cao}. As its applications, we establish a kind of hybrid-type finite continued fraction approximations related to BBP-type series of the constant $\pi$ and other classical constants, such as Catalan constant, $\pi^2$, etc.
\end{abstract}

\section{Introduction}
In the theory of mathematical constants(for example, $\pi$,  Euler-Mascheroni constant $\gamma$,  Catalan constant $G$, $\ln 2$, etc.), it is very important to construct new sequences which converge to these fundamental constants with increasingly high speed. See e.g. the survey paper of Bailey, Borwein, Mattingly, and Wightwick~\cite{BBMW} and references therein, and the books of Brent and Zimmermann~\cite{BZ}, Graham, Knuth and Patashnik~\cite{GKP}, Ifrah~\cite{Ifr}, and Wilf~\cite{Wilf}. In a celebrated paper of Bailey-Borwein-Plouffe~\cite{BBP}, they proposed the following fast series
 \begin{align}
\pi=\sum_{m=0}^{\infty}
\frac{1}{16^m}\left(\frac{4}{8m+1}-\frac{2}{8m+4}-\frac{1}{8m+5}
-\frac{1}{8m+6}\right)
:=\sum_{m=0}^{\infty}
\frac{\rho(m)}{16^m}.\label{pi-BBP}
\end{align}
This formula has the remarkable property that permits one to directly calculate binary digits of $\pi$, beginning at an arbitrary position $d$, without needing to calculate any of the first $d-1$
digits. Since this discovery in 1997, many BBP-type formulas for various mathematical constants have been discovered with the general form
\begin{align}
\alpha=\sum_{m=0}^{\infty}\frac{1}{b^m}\frac{p(m)}{q(m)},
\end{align}
where $\alpha$ is the constant,  $p$ and $q$ are polynomials in integer coefficients, and $b\ge 2$ is an integer numerical base.

Motivated by the important work of Mortici~\cite{Mor-1}, in this paper we will continue our previous works~\cite{CXY,Cao}, and apply the \emph{multiple-correction method} to construct some new sequences from a BBP-type series, which have faster rate of convergence. We give some examples to illustrate this method, such as $\pi$, Catalan constant $G$, $\pi^2$, etc. Moreover, we establish sharp bounds for the related error terms. It should be stressed that the investigation of the error terms in the approximations generated by BBP-type series is very important topic, because these error estimates can be used to study
the irrationality, transcendentally involved the constants. For example, $e$(for more details, see Aigner and Ziegler~\cite{AZ}), Apr\'ey constant $\zeta(3)$(see~\cite{Ap}), etc.

The paper is organized as follows. In Section 2, we explain how to find a finite continued fraction approximation by using the
\emph{multiple-correction method}. In Section 3, Section 4 and Section 5, we discuss $\pi$, Catalan constant and $\pi^2$, respectively. In the last section, we give three further results, and analyze the related perspective of research in this direction.

\noindent {\bf Notation.} Throughout the paper, the notation
$P_k(x)$(or $Q_k(x)$) denotes a polynomial of degree $k$
in $x$. The notation $\Psi(k;x)$
means a polynomial of degree $k$ in $x$ with all of its non-zero coefficients positive, which may be different at each occurrence.
While, $\Phi(k;x)$ denotes a polynomial of degree $k$ in $x$ with
the leading coefficient equals one, which may be different at
different section. To save space, we also use the shorthand notation to write a continued fraction
\begin{align}
\frac{a_0}{b_0+}\frac{a_1}{b_1+}\frac{a_2}{b_2+}
\cdots \frac{a_K}{b_K}
=\K_{k=0}^{K}
\frac{a_k}{b_k}.
\end{align}

\section{The multiple-correction method}

Let a series  $\sum_{m=0}^{\infty}t_m$ converge to constant $\alpha$. If we use the finite sum $\sum_{m=0}^{n-1}t_m$ to approximate or compute constant $\alpha$ for some ``comparative large" positive integer $n$, the error term $E(n)$ equals to $\sum_{m=n}^{\infty}t_m$. To evaluate it more accurately, in general, we need to ``separate" extra main-term $MC(n)$ from $E(n)$ such that the new error term $E(n)-MC(n)$ has a faster rate of convergence than $E(n)$ when $n$ tends to infinity. The idea of the \emph{multiple-correction method} is that we can achieve it in some cases by looking for
the proper structure of $MC(n)$ , where $MC(n)$ is a finite
continued fraction(see~\cite{Cao}) or a \emph{Hyper-power
expansion}(see~\cite{CXY}) in $n$. Hence, in some senses, we can view it as a rational function approximation problem of the error term $E(n)$. In fact, the
\emph{multiple-correction method} is a recursive algorithm, and one of its advantages is that by repeating correction process we always can accelerate the convergence. To describe this method clearly, we will give some definitions as follows.

\bigskip
\noindent {\bf Definition 1.} We call the integer $l-m$ to be the \emph{degree} of a rational function $R(k)=\frac{P_l(k)}{Q_m(k)}$ in $k$, and write $\deg R(k)=l-m$.

\noindent {\bf Definition 2.} Let a series $\sum_{m=0}^{\infty}t_m$ be convergent. A function $t_m$ is said to be a \emph{proper BBP-type term} if it can
be written in the form
\begin{align}
t_m=R(m)\frac{\prod_{i=1}^{uu}(a_im+c_i)!}
{\prod_{j=1}^{vv}(b_jm+d_i)!}
\frac{1}{q^m},\label{tm-def}
\end{align}
in which $q\in (0,+\infty)$ is a specific constant, and

1. $R(m)$ is a rational function in $m$,

2. the $a_i,c_i, b_j$ and $d_j$ are specific integers with
$a_i>0, b_j>0$, and

3. the quantities $uu$ and $vv$ are finite, nonnegative, specific integers.

\bigskip
Throughout the paper, we always assume that the $t_m$ is
a \emph{proper BBP-type term} and $q\neq 1$.

Now we can describe the \emph{multiple-correction method} as the following steps:
\bigskip

(Step 1) Simplify the ratio $\frac{t_{m+1}}{t_m }$ to bring the form $\frac{P_r(m)}{Q_s(m)}$, where $P, Q$ are polynomials.

(Step 2) We begin from $k=0$, and  in turn find  the finite continued approximation solution $MC_k(m)$ of the difference equation
\begin{align}
y(m)-\frac{P_r(m)}{Q_s(m)}y(m+1)-R(m)=0,\label{First-DF}
\end{align}
until some suitable $k=k^*$ you want.

(Step 3) Substitute the above $k$-th correction function $MC_k(m)$ into the left-hand side of~\eqref{First-DF} to find the constant $C_k$ and positive integer $K_0$ such that
\begin{align}
MC_k(m)-\frac{P_r(m)}{Q_s(m)}MC_k(m+1)-R(m)
+\frac{C_k}{m^{K_0}}=O\left(\frac{1}{m^{K_0+1}}\right),
\label{Second-DF}
\end{align}

(Step 4) Consider the new \emph{proper BBP-type term} appearing in (Step 3)
\begin{align}
tt_m=\frac{1}{m^{K_0}}\frac{\prod_{i=1}^{uu}(a_im)!}
{\prod_{j=1}^{vv}(b_jm)!}\frac{1}{q^m},\label{ttm-def}
\end{align}
then repeat (Step 1) to (Step 3), here it should be noted that it is often suffice for us to obtain some weak results in these cases.

(Step 5) Define the $k$-th correction error term $E_k(n)$ as
\begin{align}
E_k(n):=\alpha-\sum_{m=0}^{n-1}t_m
-\frac{\prod_{i=1}^{uu}(a_in+c_i)!}
{\prod_{j=1}^{vv}(b_jn+d_i)!}
\frac{1}{q^n}\mathrm{MC}_k(n).
\end{align}

Prove the rate of convergence of the $k$-th correction error term $E_k(n)$ when $n$ tends to infinity.

(Step 6) Based on (Step 5), we further prove sharp double-sides inequalities of $E_k(n)$ for as possible as smaller $n$.

\bigskip

Here it should be worth remarking that (Step 2) plays an important role in the \emph{multiple-correction method}. The idea of the above algorithm is originated from Mortici~\cite{Mor-1} and Gosper's Algorithm(see Chapter 5 of ~ Petkovsek, Wilf and  Zeilberger\cite{PWZ}).

Now we explain how to look for all the related coefficients in $MC_k(m)$. The initial-correction function $MC_0(m)$ is
vital. Let $\deg MC_0(m)=-\kappa_0\in \mathbb{Z}$, and denote its first coefficient by $\lambda_0\neq 0$. It is not difficult to obtain $\kappa_0$ and $\lambda_0$, which satisfy the following condition:
\begin{align}
\min_{\kappa,\lambda}\deg\left(\frac{\lambda}{m^{\kappa}}
-\frac{P_r(m)}{Q_s(m)}\frac{\lambda}{(m+1)^{\kappa}}-R(m)\right).
\end{align}

If $\kappa_0>0$, then $MC_0(m)$ has the form $\frac{\lambda_0}{\Phi(\kappa_0;m)}$. Otherwise, it is a polynomial
of degree $-\kappa_0$ with the leading coefficient $\lambda_0$. Next, just did as our previous paper~\cite{CXY,Cao}(also see (2.7) below), we can  find other coefficients in $MC_0(m)$ by solving a linear equation in turn.

Once one determines the initial-correction function $MC_0(m)$, other correction functions $MC_k(m)$ for $k\ge 1$ will become easy. Actually, one may apply two approaches to treat them. One method is power
series expansion, another is that putting the whole thing over a common denominator such that
\begin{align}
\deg\left(MC_k(m)-\frac{P_r(m)}{Q_s(m)}MC_k(m+1)-R(m)\right)
\end{align}
is a strictly decreasing function of $k$.

Next, we explain how to do (Step 5). First, by multiplying the formula (2.3) by $\frac{\prod_{i=1}^{uu}(a_im+c_i)!}
{\prod_{j=1}^{vv}(b_jm+d_i)!}
\frac{1}{q^m}$, then by adding these formulas from $m=n$ to $m=\infty$, finally by checking $$\lim_{n\rightarrow\infty}\frac{\prod_{i=1}^{uu}(a_in+c_i)!}
{\prod_{j=1}^{vv}(b_jn+d_i)!}
\frac{1}{q^n}\mathrm{MC}_k(n)=0,$$
in this way it is not difficult to get the desired results for the rate of convergence of the $k$-th correction error term $E_k(n)$.

Finally, there doesn't exist the general method to treat (Step 6), which needs many delicate estimations for the involved series.

Since $MC_k(m)$ and other constants need a huge of computations, we often use an appropriate symbolic computation software. In addition, the exact expression at each occurrence also takes a lot of space. Hence, in this paper we omit some related details for space limitation. For interested readers can see our previous papers~\cite{CXY,Cao,XY}.
\begin{rem}
If $q=1$, and $uu=vv=0$, we may replace equation~\eqref{First-DF} by $y(m)-y(m+1)-R(m)=0$, the above method is still efficient.
\end{rem}

\noindent {\bf An example.}  We would like to give an example to show how to manipulate (Step 1) to (Step 3). It is well-known that
\begin{align}
\frac{1}{\pi}=\frac{1}{16}\sum_{m=0}^{\infty}
\frac{((2m)!)^3}{(m!)^6}\frac{42m+5}{4096^m}:
=\frac{1}{16}\sum_{m=0}^{\infty}t_m,
\end{align}
which is proposed by Srinivasa Ramanujan~\cite{Ra}, also see (1.4) of Mortici~\cite{Mor-1}.
We take $R(m)=42m+5, q=4096, uu=3$ and $vv=6$ in Definition 2, hence it is a \emph{proper BBP-type term}.

(Step 1) It is easy to check
\begin{align}
\frac{t_{m+1}}{t_m}=\frac{(2m+2)^3(2m+1)^3}{4096(m+1)^6}.
\end{align}

(Step 2) We choose $k^*=6$. Consider the difference equation
\begin{align}
y(m)-\frac{(2m+2)^3(2m+1)^3}{4096(m+1)^6}y(m+1)-R(m)=0.
\end{align}
By using \emph{Mathematica} software, it is not difficult to find
\begin{align}
\mathrm{MC}_0(m)=\frac{128}{3}m+\frac{128}{27},\quad
\mathrm{MC}_k(m)=\mathrm{MC}_0(m)+
\K_{j=1}^{k}\frac{a_j}{m+b_j},(k\ge 1),
\end{align}
where
\begin{align*}
a_1=&\frac{32}{81}, &b_1=\frac{10}{9},\\
a_2=&-\frac{7}{324}, &b_2=\frac{27}{7},\\
a_3=&\frac{19856}{3969}, &b_3=-\frac{145795}{156366},\\
a_4=&\frac{1396171}{4620243}, &b_4=\frac{15549372115}{4455381114},\\
a_5=&-\frac{818973874600}{3222301435929}, &b_5=\frac{24496617933181}{3948754138854},\\
a_6=&\frac{7676419604757068}{881904503553129},
&b_6=-\frac{535521415681420831}{571477212182467206}.
\end{align*}
(Step 3) By using \emph{Mathematica} software again, we easily find
$K_0=2k+1$ and
$
C_0=\frac{7}{18},
C_1=\frac{49}{5832},
C_2=-\frac{2482}{59049},
C_3=\frac{2792342}{219839427},
C_4=\frac{9239028400}{2861932547871},
C_5=-\frac{15394944382400}{547865648052993},
C_6=\frac{5377668984891011200}{100647847362777935517},
$
which satisfy
\begin{align}
\mathrm{MC}_k(m)-\frac{(2m+2)^3(2m+1)^3}{4096(m+1)^6}
\mathrm{MC}_k(m+1)-R(m)+\frac{C_k}{m^{2k+1}}=O\left(
\frac{1}{m^{2k+2}}
\right).
\end{align}
Now we let
\begin{align}
E_k(n):=\frac{1}{\pi}-\frac{1}{16}\sum_{m=0}^{n-1}
\frac{((2m)!)^3}{(m!)^6}\frac{42m+5}{4096^m}
-\frac{((2n)!)^3}{16(n!)^6}\frac{\mathrm{MC}_k(n)}
{4096^{n}}.\label{1/pi-Ek-Def}
\end{align}
Then for $0\le k\le 6$, we may prove
\begin{align}
\lim_{n\rightarrow\infty}\frac{(n!)^6}{((2n)!)^3}
4096^nn^{2k+1}E_k(n):=\frac{4C_k}{63}.
\end{align}

\section{The results for $\pi$}
In order to illustrate the so-called the \emph{multiple-correcton method} formulated in previous section, first we will prove the following theorem.
\begin{thm} Let $\rho(m)$ be defined as ~\eqref{pi-BBP}. For every integer $k\ge 0$, the $k$-th correction $\mathrm{MC}_k(n)$ is defined by
\begin{align}
\mathrm{MC}_0(n)=&\frac{\frac{1}{4}}{(n+\frac{7}{16})^2
-\frac{73}{256}},\quad \mathrm{MC}_k(n)=\frac{\frac 14}{(n+\frac{7}{16})^2-\frac{73}{256}+}
\K_{j=1}^{k}\frac{a_j}{n+b_j}
, (k\ge 1),
\end{align}
where
\begin{align*}
a_1=&\frac{21}{64},&b_1=\frac{15}{7},\\
a_2=&-\frac{265}{392},&b_2=\frac{9299}{2968},\\
a_3=&-\frac{3381}{2809},&b_3=\frac{20517}{4876},\\
a_4=&-\frac{18921}{8464},&b_4=\frac{94519}{21896},\\
a_5=&-\frac{3260043}{453152},&b_5=\frac{25408967}{7496524},\\
a_6=&-\frac{3740382415}{496062002},&b_6=\frac{482484243355}
{72002790104},\\
a_7=&-\frac{435259601465}{326597391169},
&b_7=\frac{133863589556959}{4859799720860},\\
a_8=&\frac{18170745077870217}{36157137144200},
&b_8=-\frac{550189873911066313}{30042487323672220},\\
a_9=&\frac{1184188272901493239625}{399390489791710771232},
&b_9=\frac{55409761792537711960915}{5291704918098810592904}.
\end{align*}
Let the $k$-th correction error term $E_k(n)$ be defined as
\begin{align}
E_k(n):=\pi -\sum_{m=0}^{n-1}
\frac{\rho(m)}{16^m}-\frac{1}{16^n}\mathrm{MC}_k(n),
\end{align}
Then for all integers $0\le k\le 9$, we have
\begin{align}
\lim_{n\rightarrow\infty}16^{n}n^{2k+5}E_k(n)=\frac{16C_k}{15},
\end{align}
where
$
C_0=-\frac{315}{4096},
C_1=-\frac{11925}{229376},
C_2=-\frac{108675}{1736704},
C_3=-\frac{1686825}{12058624},
C_4=-\frac{287025525}{285212672},
C_5=-\frac{4009909971375}{528448749568},
C_6=-\frac{27702923551875}{2739417382912},
C_7=\frac{580053423565590975}{114135803101184},
C_8=-\frac{457280686810171702603125}{30346857643463671808},
C_9=-\frac{197080602286603349404715625}{1608316872287169019904}.
$
\end{thm}

\begin{lem}  Under the same notation of Theorem 1, when $m$ tends to $\infty$, we have for $0\le k\le 9$
\begin{align}
\mathrm{MC}_k(m)-\frac{1}{16}\mathrm{MC}_k(m+1)-\rho(m)
+\frac{C_{k}}{m^{2k+5}}=O\left(\frac{1}{m^{2k+6}}\right).
\label{pi-difference-inequality}
\end{align}
\end{lem}
\proof By using \emph{Mathematica} software, we expand $\mathrm{MC}_k(m)-\frac{1}{16}\mathrm{MC}_k(m+1)-\rho(m)$ as the power series in terms of $m^{-1}$, then after some simplifications we can prove \eqref{pi-difference-inequality}.
\begin{lem} Let
\begin{align}
u(n)=\frac{1}{\frac{15}{16}n^{23}+\frac{23}{16}n^{22}},\quad v(n)=\frac{1}{\frac{15}{16}n^{23}+\frac{23}{16}n^{22}- \frac{4163}{240}n^{21}}.
\end{align}
Then for $n\ge 4$, we have
\begin{align}
\frac{1}{16^n}u(n)<\sum_{m=n}^{\infty}\frac{1}{m^{23}16^m}
<\frac{1}{16^n}v(n).\label{pi-l-2}
\end{align}
\end{lem}
\proof We can check for $m\ge 4$
\begin{align}
&u(m)-\frac{1}{16}u(m+1)-\frac{1}{m^{23}}=
-\frac{\Psi_1(22;m)}{m^{23} (23+15m)\Psi_2(23;m)}<0,\\
&v(m)-\frac{1}{16}v(m+1)-\frac{1}{m^{23}}\\
&=\frac{
\Psi_{3}(21;m)(m-4)+4593130341153628118305}{m^{23}((225 m+1245)(m-4)817)(\Psi_4(22;m)(m-4)+1519680023193359375)}>0.
\nonumber
\end{align}
By multiplying (3.7) and (3.8) by $16^{-m}$, we obtain the telescoping inequalities
\begin{align*}
\frac{1}{16^{m}}u(m)-\frac{1}{16^{m+1}}u(m+1)
<\frac{1}{m^{23}16^m}<\frac{1}{16^{m}}v(m)-\frac{1}{16^{m+1}}v(m+1).
\end{align*}
Now by adding the above inequalities from $m=n$ to $m=\infty$, we can obtain \eqref{pi-l-2} at once.

\noindent{{\emph{The proof of Theorem 1.}}
First, by multiplying \eqref{pi-difference-inequality} by $16^{-m}$, we have
\begin{align}
\frac{1}{16^m}\mathrm{MC}_k(m)-\frac{1}{16^{m+1}}\mathrm{MC}_k(m+1)
-\frac{\rho(m)}{16^m}+\frac{C_k}{16^mm^{2k+5}}=
O\left(\frac{1}{16^mm^{2k+6}}\right).
\end{align}
Then, by adding these formulas from $m=n$ to $m=\infty$, we get
\begin{align}
E_k(n)=&\sum_{m=n}^{\infty}\frac{\rho(m)}{16^m}-
\frac{\mathrm{MC}_k(n)}{16^n}\\
=&\sum_{m=n}^{\infty}\frac{C_k}{16^mm^{2k+5}}
+O\left(\sum_{m=n}^{\infty}\frac{1}{16^mm^{2k+6}}\right).\nonumber
\end{align}
It is easy to prove
\begin{align}
\sum_{m=n}^{\infty}\frac{1}{16^mm^{2k+6}}
=O\left(\frac{1}{n^{2k+6}}\sum_{m=n}^{\infty}\frac{1}{16^m}\right)
=O\left(\frac{1}{n^{2k+6}16^n}\right).
\end{align}
Combining (3.10), (3.11) and \eqref{pi-l-2} completes the the proof of Theorem 1 in case of $k=9$. For $0\le k\le 8$, we may prove the theorem in the same approach.

The following theorem tells us how to improve (3.3).
\begin{thm} Under the same notation of Theorem 3, we have for $n\ge 88$
\begin{align}
\frac{16C_9}{15\cdot 16^n(n+1)^{23}}<E_9(n)<\frac{16C_9}{15\cdot 16^n(n+5)^{23}}.
\end{align}
\end{thm}

\begin{lem}
Let $f(m)=\mathrm{MC}_9(m)-\frac{1}{16}\mathrm{MC}_9(m+1)-\rho(m)$. Then we have for $m\ge 41$
\begin{align}
-\frac{D_{10}}{(m+\frac {55}{32})^{24}}<f(m)+\frac{C_9}{m^{23}}<
-\frac{D_{10}}{(m+\frac {71}{32})^{24}},\label{Catalan-lemma 3}
\end{align}
where
$D_{10}=\frac{28928763399211176287296777111194638413125}
{2409036421853659622126333496131584}$.
\end{lem}
\proof We can check for $m\ge 41$
\begin{align}
&f(m)+\frac{C_9}{m^{23}}+\frac{D_{10}}{(m+\frac{55}{32})^{24}}\\
&=\frac{\Psi_5(47;m)(m-41)+702881\cdots 544759}{301\cdots 448m^{23}(1+2m)(3+4m)(1+8m)
(5+8m)(55+32m)^{24}\Psi_6(22;m)}>0.\nonumber
\end{align}
Similarly, we can check for $m\ge 1$
\begin{align}
&f(m)+\frac{C_9}{m^{23}}+\frac{D_{10}}{(m+\frac {71}{32})^{24}}\\
&=-\frac{162820783125\Psi_7(48;m)}{301\cdots 448m^{23}(1+2m)(3+4m)(1+8m)
(5+8m)(71+32m)^{24}\Psi_8(22;m)}<0.\nonumber
\end{align}
This completes the proof of Lemma 3.
\begin{lem} We let
\begin{align}
&u_1(n)=\frac{1}{\frac{15}{16}(n+\frac{63}{32})^{24}},\quad v_1(n)=\frac{1}{\frac{15}{16}(n+\frac{55}{32})^{24}},\\
&u_2(n)=\frac{1}{\frac{15}{16}(n+\frac{79}{32})^{24}},\quad
v_2(n)=\frac{1}{\frac{15}{16}(n+\frac{71}{32})^{24}}.
\end{align}
Then for all positive integer $n$
\begin{align}
&\frac{1}{16^n}u_1(n)<\sum_{m=n}^{\infty}\frac{1}{(m+\frac {55}{32})^{24}16^{m}}
<\frac{1}{16^n}v_1(n),
\\
&\frac{1}{16^n}u_2(n)<\sum_{m=n}^{\infty}\frac{1}{(m+\frac {71}{32})^{24}16^{m}}
<\frac{1}{16^n}v_2(n).
\end{align}
\end{lem}
\proof Similar to the proof of (3.11), we can prove the inequalities of right-hand sides in both (3.18) and (3.19) trivially. By using \emph{Mathematica} software, it is not difficult to prove
\begin{align}
&u_1(m)-\frac{1}{16}u_1(m+1)-\frac{1}{(m+\frac {55}{32})^{24}}\\
=&-\frac{8507\cdots2864\Psi_1(47;m)}
{15(55+32m)^{24}(63+32m)^{24}(95+32m)^{24}}<0,\nonumber\\
&u_2(m)-\frac{1}{16}u_2(m+1)-\frac{1}{(m+\frac {71}{32})^{24}}\\
=&-\frac{8507\cdots2864\Psi_2(47;m)}
{15(71+32m)^{24}(79+32m)^{24}(111+32m)^{24}}<0,\nonumber
\end{align}
Now  by multiplying  the above two inequalities by $16^{-m}$, then by adding these formulas from $m=n$ to $m=\infty$, we can prove the other two inequalities in Lemma 4 immediately.

\noindent{{\emph{The proof of Theorem 2.}} By multiplying  (3.13) by $16^{-m}$ , then by adding these formulas from $m=n$ to $m=\infty$, we have
\begin{align}
\sum_{m=n}^{\infty}\frac{C_9}{16^mm^{23}}
+\sum_{m=n}^{\infty}\frac{D_{10}}{16^m(m+\frac{71}{32})^{24}}
<E_9(n)<\sum_{m=n}^{\infty}\frac{C_9}{16^mm^{23}}
+\sum_{m=n}^{\infty}\frac{D_{10}}{16^m(m+\frac{55}{32})^{24}}.
\end{align}
By using \emph{Mathematica} software, it is not difficult to prove for $n\ge 88$
\begin{align}
&\frac{1}{15n^{23}+23n^{22}}
+\frac{D_{10}}{C_9}\frac{1}{15(n+\frac{47}{32})^{24}}-
\frac{1}{15(n+5)^{23}}\\
&=\frac{\Psi_3(45;n)(n-88)+147532\cdots 759049}{3399\cdots 8325 n^{22}(5+n)^{23}(55+32n)^{24}(23+15n)}>0.\nonumber
\end{align}
By (3.23), Lemma 2 and Lemma 4, we have
\begin{align}
&\sum_{m=n}^{\infty}\frac{C_9}{16^mm^{23}}
+\sum_{m=n}^{\infty}\frac{D_{10}}{16^m(m+\frac{55}{32})^{24}}\\
<& \frac{C_9}{16^n}u(n)+\frac{D_{10}}{16^n}v_1(n)\nonumber\\
=&\frac{16C_9}{16^n}\left(\frac{1}{15n^{23}+23n^{22}}
+\frac{D_{10}}{C_9}\frac{1}{15(n+\frac{55}{32})^{24}}\right)
\nonumber\\
<&\frac{16C_9}{15\cdot 16^n(n+5)^{23}},\quad (n\ge 88).\nonumber
\end{align}

Similarly, we also have from Lemma 2 and Lemma 4
\begin{align}
&\sum_{m=n}^{\infty}\frac{C_9}{16^mm^{23}}
+\sum_{m=n}^{\infty}\frac{D_{10}}{16^m(m+\frac{71}{32})^{24}}\\
>&\frac{C_9}{16^n}v(n)+\frac{D_{10}}{16^n}u_2(n)\nonumber\\
=&\frac{16C_9}{16^n}\left(\frac{1}{15n^{23}+23n^{22}- \frac{4163}{15}n^{21}}
+\frac{D_{10}}{C_9}\frac{1}{15(n+\frac{79}{32})^{24}}\right)\nonumber
>\frac{16C_9}{15\cdot 16^n(n+1)^{23}}
\end{align}
Here we note for $n\ge 41$
\begin{align}
&\frac{1}{15n^{23}+23n^{22}- \frac{4163}{15}n^{21}}
+\frac{D_{10}}{C_9}\frac{1}{15(n+\frac{79}{32})^{24}}-
\frac{1}{15(n+1)^{23}}\\
&=-\frac{\Psi_4(45;n)(n-41)+58459\cdots 68269}{3399\cdots 8325(1+n)^{23}(79+32n)^{24}n^{21}\left((225 n+9570)(n-41)+388207\right)}<0.\nonumber
\end{align}
Finally, Theorem 2 follows from (3.22), (3.24) and (3.25) at once.

\section{Catalan constant}
Catalan constant can be defined as
\begin{align}
G=\sum_{m=0}^{\infty}\frac{(-1)^m}{(2m+1)^2}=0.915965594\cdots,
\end{align}
which is arguably the most basic constant whose irrationality and transcendence remain unproven. The most economical BBP-type series for computing Catalan constant may be
\begin{align}
G=&\frac{1}{4096}\sum_{m=0}^{\infty}\frac{1}{4096^m}\left(
\frac{36864}{(24 m+2)^2}
-\frac{30720}{(24 m+3)^2}
-\frac{30720}{(24 m+4)^2}
-\frac{6144}{(24 m+6)^2}\right.\label{Catalan-BBP-formula}\\
&-\frac{1536}{(24 n+7)^2}
+\frac{2304}{(24 m+9)^2}
+\frac{2304}{(24 m+10)^2}
+\frac{768}{(24 m+14)^2}\nonumber\\
&+\frac{480}{(24 m+15)^2}
+\frac{384}{(24 m+11)^2}
+\frac{1536}{(24 m+12)^2}
+\frac{24}{(24 m+19)^2}\nonumber\\
&\left.-\frac{120}{(24 m+20)^2}
-\frac{36}{(24 m+21)^2}
+\frac{48}{(24 m+22)^2}
-\frac{6}{(24 m+23)^2}\right).\nonumber\\
:=&\frac{1}{4096}\sum_{m=0}^{\infty}\frac{\nu(m)}{4096^m},\nonumber
\end{align}
see formula (18) in Bailey, Borwein, Mattingly and Wightwick~\cite{BBMW}.

\begin{thm}
For every  integer $k\ge 0$, the $k$-th correction
function $\mathrm{MC}_k(n)$ is defined by
\begin{align}
\mathrm{MC}_0(n)=\frac{-\frac{3}{128}}{(n + \frac{13}{72})^2 + \frac{41}{432}},\quad
\mathrm{MC}_k(n)=\frac{-\frac{3}{128}}{(n + \frac{13}{72})^2 + \frac{41}{432}+}\K_{j=1}^{k}\frac{a_j}{n+b_j},(k\ge 1),
\end{align}
where
\begin{align*}
a_1=&-\frac{517}{23328}, &b_1=\frac{156655}{148896},\\
a_2=&\frac{366823315}{821111808}, &b_2=-\frac{73939238279831}{163855572930720},\\
a_3=&-\frac{975884794104398189}{98093762087712545025}, &b_3=-\frac{1932406340618716298628867667}
{299122227350085279497481360},\\
a_4=&\frac{1518828040567790867982188908085299115}
{24627466973909279332577879325543168},\\ b_4=&\frac{10414320422851149518238529301402392329619615007}
{1125439555781241535752796061860324225756510608}.
\end{align*}
We define the $k$-th correction error term $E_k(n)$ as
\begin{align}
E_k(n):=G-\frac{1}{4096}\sum_{m=0}^{n-1}\frac{\nu(m)}{4096^m}
-\frac{1}{4096^{n+1}}\mathrm{MC}_k(n).\label{Catalan-Ek-Def}
\end{align}
Then for all integers $0\le k\le 4$, we have
\begin{align}
\lim_{n\rightarrow\infty}4096^{n}n^{2k+5}E_k(n)=\frac{C_k}{4095},
\end{align}
where
$
C_0=-\frac{235235}{248832},
C_1=\frac{166904608325}{395200954368} ,
C_2=-\frac{171770824494197747}{40882919041278148608},
C_3=\frac{1883922668487810936804537501055}{
7270540656226904507330240446464},
C_4=-\frac{59816694319657990230589749754634406261775}{
191377827680729835340592443669705291028496384}.
$
\end{thm}
\begin{lem} Let $0\le k\le 4$. Under the same notation of Theorem 3, when $m$ tends to infinity we have
\begin{align}
\mathrm{MC}_k(m)-\frac{1}{4096}\mathrm{MC}_k(m+1)-\nu(m)
+\frac{C_k}{m^{2k+5}}=O\left(\frac{1}{m^{2k+6}}\right).
\label{Catalan-difference-inequality}
\end{align}
\end{lem}
\proof First, by using \emph{Mathematica} software we expand $\mathrm{MC}_k(m)-\frac{1}{4096}\mathrm{MC}_k(m+1)-\nu(m)$ as the power series in terms of $m^{-1}$, then after some simplifications we can prove \eqref{Catalan-difference-inequality}.

\bigskip
\begin{lem}
We let
\begin{align}
u(n):=\frac{1}{\frac{4095}{4096}n^{13} + \frac{13}{4096}n^{12}},\quad
v(n):=\frac{1}{\frac{4095}{4096}n^{13} + \frac{13}{4096}n^{12}-\frac{14333}{645120}n^{11}}.
\end{align}
Then for all positive integers $n$, we have
\begin{align}
\frac{1}{4096^n}u(n)<\sum_{m=n}^{\infty}\frac{1}{m^{13}4096^m}
<\frac{1}{4096^n}v(n).
\end{align}
\end{lem}
\proof  By applying \emph{Mathematica} software, we can prove for $m\ge 1$
\begin{align}
u(m)-\frac{1}{4096}u(m+1)-\frac{1}{m^{13}}=
-\frac{\Psi_1(12;m)}{13m^{13}(1+m)^{12}(1+315m)(316+315m)}<0,
\label{u(n)-inequality}
\end{align}
and
\begin{align}
&v(m)-\frac{1}{4096}v(m+1)-\frac{1}{m^{13}}\label{v(n)-inequality}\\
=&
\frac{\Psi_2(12;m)}{m^{13}(1+m)^{11}(-28666+4095m+1289925m^2)
(1265354+2583945m+1289925m^2)}>0.\nonumber
\end{align}
We multiply \eqref{u(n)-inequality} by $\frac{1}{4096^m}$ to obtain the telescoping inequality
\begin{align*}
\frac{1}{4096^m}u(m)-\frac{1}{4096^{m+1}}u(m+1)
-\frac{1}{m^{13}4096^m}<0.
\end{align*}
Then, by adding these inequalities from $m=n$ to $m=\infty$, we have for all integers $n\ge 1$
\begin{align}
\frac{1}{4096^n}u(n)-\sum_{m=n}^{\infty}\frac{1}{m^{13}4096^m}
<0.\label{lemma2-lower bound}
\end{align}
Similarly, we multiply \eqref{v(n)-inequality} by $\frac{1}{4096^m}$ to get the telescoping inequality
\begin{align*}
\frac{1}{4096^m}v(m)-\frac{1}{4096^{m+1}}v(m+1)
-\frac{1}{m^{13}4096^m}>0.
\end{align*}
Then, by adding these inequalities from $m=n$ to $m=\infty$, we have for all integer $n\ge 1$
\begin{align}
\frac{1}{4096^n}v(n)-\sum_{m=n}^{\infty}\frac{1}{m^{13}4096^m}
>0.\label{lemma2-upper bound}
\end{align}
Finally, combining ~\eqref{lemma2-lower bound} and ~\eqref{lemma2-upper bound} completes the proof of Lemma 6.

\bigskip
\noindent {{\emph{The proof of Theorem 3}}. We only give the proof of Theorem 1 in the case of $k=4$, the other can be proved similarly. By multiplying \eqref{Catalan-difference-inequality} by $\frac{1}{4096^m}$, we
get the telescoping estimate
\begin{align*}
\frac{1}{4096^m}\mathrm{MC}_4(m)-\frac{1}{4096^{m+1}}
\mathrm{MC}_4(m+1)
-\frac{\nu(m)}{4096^m}+\frac{C_4}{4096^mm^{13}}
=O\left(\frac{1}{4096^mm^{14}}\right).
\end{align*}
Then, by adding these formulas from $m=n$ to $m=\infty$, we have
\begin{align}
\frac{1}{4096^n}\mathrm{MC}_4(n)-\sum_{m=n}^{\infty}
\frac{\nu(m)}{4096^m}
+\sum_{m=n}^{\infty}\frac{C_4}{4096^mm^{13}}=
O\left(\frac{1}{n^{14}}\sum_{m=n}^{\infty}\frac{1}{4096^m}\right).
\end{align}
It is not difficult to check that
\begin{align*}
\sum_{m=n}^{\infty}\frac{1}{m^{14}4096^m}=O\left(\frac{1}{n^{14}}
\sum_{m=n}^{\infty}\frac{1}{4096^m}\right)
=O\left(\frac{1}{n^{14}4096^n}\right).
\end{align*}
Thus from (4.13)
\begin{align}
\sum_{m=n}^{\infty}\frac{\nu(m)}{4096^m}-\frac{1}{4096^n}
\mathrm{MC}_4(n)
=&\sum_{m=n}^{\infty}\frac{C_4}{4096^mm^{13}}
+O\left(\frac{1}{4096^nn^{14}}\right)
\label{Catalan-error-estimate}\\
=&C_4\sum_{m=n}^{\infty}\frac{1}{4096^mm^{13}}
+O\left(\frac{1}{4096^nn^{14}}\right).
\nonumber
\end{align}
It follows from ~\eqref{Catalan-BBP-formula},~\eqref{Catalan-Ek-Def},
\eqref{Catalan-error-estimate}
\begin{align}
E_4(n)=&\frac{1}{4096}\sum_{m=n}^{\infty}\frac{\nu(m)}{4096^m}
-\frac{1}{4096^{n+1}}\mathrm{MC}_4(n)\label{Catalan-Ek-estimate}\\
&=\frac{C_4}{4096}\sum_{m=n}^{\infty}\frac{1}{4096^mm^{13}}
+O\left(\frac{1}{4096^nn^{14}}\right).\nonumber
\end{align}
Now combining~\eqref{Catalan-Ek-estimate} and Lemma 6 finishes the
proof of Theorem 3.

\begin{thm} Under the same notation of Theorem 3,
we have the following double-sides inequalities for $n\ge 12$
\begin{align}
\frac{C_4}{4095\cdot 4096^nn^{13}}<E_4(n)<\frac{C_4}{4095\cdot 4096^n(n+5)^{13}}.\label{Catalan-Ek-inequalities}
\end{align}
\end{thm}
\begin{rem}
In fact, by applying the same method as the proof Theorem 4, in the cases of $0\le k\le 3$ we can get analogous estimates for $E_k(n)$ . Here we leave these for readers to check.
\end{rem}

\begin{lem}
Let $f(n)=\mathrm{MC}_4(m)-\frac{1}{4096}\mathrm{MC}_4(m+1)-\nu(m)$. Then we have for $m\ge 2$
\begin{align}
-\frac{D_5}{(m+\frac 14)^{14}}<f(m)+\frac{C_4}{m^{13}}<
-\frac{D_5}{(m+\frac 34)^{14}},\label{Catalan-lemma 3}
\end{align}
where
$D_5=\frac{168833398318043946755965656344277667293765672161072072
7013667676284669255}{3423158769571652580795886387022473556501828
66279601922774818539625196290048}$.
\end{lem}

\proof By using \emph{Mathematica} software, we can prove for $m>0$
\begin{align}
f(m)+\frac{C_4}{m^{13}}+\frac{D_5}{(m+\frac 34)^{14}}=
-\frac{5\Psi_1(56;m)}{21394742\cdots74768128\Psi_2(71;m)}<0.
\end{align}
This completes the proof of right-hand side inequality of
\eqref{Catalan-lemma 3}.
Similarly, one has
\begin{align}
f(m)+\frac{C_4}{m^{13}}+\frac{D_5}{(m+\frac 14)^{14}}=
\frac{5\Psi_3(53;m)(m-2)+261503\cdots 734375}
{21394742\cdots 74768128m^{13}\Psi_4(56;m)}>0.
\end{align}
Hence the left-hand side inequality of
\eqref{Catalan-lemma 3} holds for $m\ge 2$. This completes the proof of Lemma 7.
\begin{lem}
For $n\ge 1$, we have
\begin{align}
&\frac{1}{4095\cdot 4096^{n-1}(n+\frac 12)^{14}}<\sum_{m=n}^{\infty}\frac{1}{4096^m(m+\frac 14)^{14}}
<\frac{1}{4095\cdot 4096^{n-1}(n+\frac 14)^{14}}.\\
&\frac{1}{4095\cdot 4096^{n-1}(n+1)^{14}}<\sum_{m=n}^{\infty}\frac{1}{4096^m(m+\frac 34)^{14}}
<\frac{1}{4095\cdot 4096^{n-1}(n+\frac 34)^{14}}.
\end{align}
\end{lem}
\proof We note that both upper bounds in the lemma are trivial.
Let
\begin{align*}
&r_1(m)=\frac{4096}{4095}\frac{1}{(m+\frac 12)^{14}},\quad
r_2(m)=\frac{4096}{4095}\frac{1}{(m+1)^{14}},\\
&s_1(m)=\frac{1}{(m+\frac 14)^{14}},\quad
s_2(m)=\frac{1}{(m+\frac 34)^{14}}.
\end{align*}

It is not difficult to check
\begin{align}
&r_1(m)-\frac{1}{4096}r_1(m+1)-s_1(m)=-
\frac{16384\Psi_5(27;m)}{4095(1+2m)^{14}(3+2m)^{14}(1+4m)^{14}}<0,\\
&r_2(m)-\frac{1}{4096}r_2(m+1)-s_2(m)=-\frac{\Psi_6(27;m)}
{4095(1+m)^{14}(2+m)^{14}(3+4m)^{14}}<0.
\end{align}
Now  by multiplying the above two inequalities by $4096^{-m}$, then by adding these formulas from $m=n$ to $m=\infty$, we can prove the other two inequalities in Lemma 8 immediately.

\noindent{{\emph{The proof of Theorem 4.}} Similar to the proof of~\eqref{Catalan-Ek-estimate}, from Lemma 7 we have
\begin{align}
\frac{D_5}{4096}\sum_{m=n}^{\infty}\frac{1}{4096^m(m+\frac 34)^{14}}
<E_4(n)-\frac{C_4}{4096}\sum_{m=n}^{\infty}\frac{1}{4096^mm^{13}}
<\frac{D_5}{4096}\sum_{m=n}^{\infty}\frac{1}{4096^m(m+\frac 14)^{14}}.
\label{Catalan-Ek-inequalities-1}
\end{align}
From Lemma 8 and Lemma 6, we have
\begin{align}
E_4(n)&<\frac{C_4}{4096}\frac{1}{4096^n}u(n)
+\frac{D_5}{4096}\frac{1}{4095\cdot 4096^{n-1}(n+\frac 14)^{14}}\label{C-E4-upper bound}\\
&=\frac{C_4}{4095\cdot 4096^n}\left(
\frac{1}{n^{13}+\frac{1}{315}n^{12}}+\frac{D_5}{C_4}
\frac{1}{(n+\frac 14)^{14}}\right)\nonumber\\
&<\frac{C_4}{4095\cdot 4096^n(n+5)^{13}}\quad (n\ge 12).\nonumber
\end{align}
Here we use $\frac{(n+5)^{13}}{n^{13}+\frac{1}{315}n^{12}}+\frac{D_5}{C_4}
\frac{(n+5)^{13}}{(n+\frac 14)^{14}}-1=
\frac{\Psi_5(13;n)(n-12)+542244\cdots 228831}{102878\cdots 080445 n^{12}(1+4n)^{14}(1+315n)}>0$ for $n\ge 12.$

Similarly, we can check that for $n>2$
\begin{align*}
&n^{13}\left(
\frac{1}{n^{13}+\frac{n^{12}}{315}-\frac{28666n^{11}}{1289925}}
+\frac{D_5}{C_4}\frac{1}{(n+1)^{14}}\right)-1\\
=&-\frac{\Psi_6(14;n)(n-2)+932326\cdots 839936}{164605\cdots 287120(1+n)^{14}(-28666+4095n+1289925n^2)}<0.
\end{align*}
Hence for $n\ge 2$
\begin{align}
E_4(n)&>\frac{C_4}{4096}\frac{1}{4096^n}v(n)
+\frac{D_5}{4096}\frac{1}{4095\cdot 4096^{n-1}(n+1)^{14}}\label{C-E4-lower bound}\\
&=\frac{C_4}{4095\cdot 4096^n}\left(
\frac{1}{n^{13}+\frac{n^{12}}{315}-\frac{28666n^{11}}{1289925}}
+\frac{D_5}{C_4}\frac{1}{(n+1)^{14}}\right).\nonumber\\
&>\frac{C_4}{4095\cdot 4096^n n^{13}}.\nonumber
\end{align}
Finally, combining~\eqref{C-E4-upper bound} and~\eqref{C-E4-lower bound} completes the proof of Theorem 4.

\section{The results for $\pi^2$ }
The following BBP-type formula is taken from (18) in Bailey, Borwein, Mattingly and Wightwick~\cite{BBMW}
\begin{align}
\pi^2=&\frac{2}{27}\sum_{m=0}^{\infty}\frac{1}{729^m}\left(
\frac{243}{(12m+1)^2}
-\frac{405}{(12m+2)^2}
-\frac{81}{(12m+4)^2}
-\frac{27}{(12m+5)^2}\right.\\
&\left.
-\frac{72}{(12m+6)^2}
-\frac{9}{(12m+7)^2}-\frac{9}{(12m+8)^2}
-\frac{5}{(12m+10)^2}
+\frac{1}{(12m+11)^2}\right)\nonumber\\
:
=&\frac{2}{27}\sum_{m=0}^{\infty}\frac{\varrho(m)}{729^m}.\nonumber
\end{align}
\begin{thm}For every integer $k\ge 0$, the $k$-th correction $\mathrm{MC}_k(n)$ is defined by
\begin{align}
\mathrm{MC}_0(n):=\frac{-\frac{10935}{5824}}
{(n+\frac{3473}{10920})^2+\frac{508433}{13249600}},
\quad
\mathrm{MC}_k(n):=\frac{-\frac{10935}{5824}}
{(n+\frac{3473}{10920})^2+\frac{508433}{13249600}+}
\K_{j=1}^{k}\frac{a_j}{n+b_j},(k\ge 1),
\end{align}
where
\begin{align*}
a_1=&\frac{1704001969}{54257112000},\quad
b_1=\frac{2133779424499}{12405134334320},\\
a_2=&-\frac{22377711469278547658588675}
{55399448826908967430750464},
\quad b_2=\frac{7838462085871364023219390913487412021}
{6662364404905290370545187619443579824},\\
a_3=&-\frac{338155884480620847387677263213133005773122041905270}
{6634895805691977782779752766105114022452560309751729},\\
b_3=&\frac{518071383229948104130947807715226040921415380062488
629146343414684409}{258632973680067531610825571620741735163688
602840992494501823795966560}.
\end{align*}
We define the $k$-th correction error term $E_k(n)$ as
\begin{align}
E_k(n):=\pi^2-\frac{2}{27}\sum_{m=0}^{n-1}\frac{\varrho(m)}{729^m}
-\frac{2}{27\cdot 729^n}\mathrm{MC}_k(n).
\end{align}
Then for $0\le k\le 3$, we have
\begin{align}
\lim_{n\rightarrow\infty}
729^nn^{2k+5}E_k(n):=\frac{27\cdot C_k}{364}.
\end{align}
where
\begin{align*}
C_0=&\frac{1704001969}{28937126400},\quad
C_1=\frac{895108458771141906343547}{37631431943365237081767936},\\
C_2=&\frac{33074676617409163665475129038532721493305}
{27282731409796850283137568847626580833927168},\\
C_3=&\frac{51782290831323026508865202336606730861855228893902257466
379094191}{256455365050612952720466033711567849946969084447258104421
96683325440}.
\end{align*}
\end{thm}
\begin{lem} Let $0\le k\le 3$.
When $m$ tends to $\infty$, we have
\begin{align}
\mathrm{MC}_k(m)-\frac{1}{729}\mathrm{MC}_k(m+1)-\varrho(m)
+\frac{C_k}{m^{2k+5}}=O\left(\frac{1}{m^{2k+6}}\right).
\end{align}
\end{lem}
\proof First, by using \emph{Mathematica} software we expand $\mathrm{MC}_k(m)-\frac{1}{729}\mathrm{MC}_k(m+1)-\varrho(m)$ as the power series in terms of $m^{-1}$, then after some simplifications we can prove Lemma 9.

\begin{lem}
Let
\begin{align}
u(n)=\frac{1}{\frac{728}{729}n^{11}+\frac{11}{729}n^{10}},\quad v(n)=\frac{1}{\frac{728}{729}n^{11}+\frac{11}{729}n^{10}-
\frac{48059}{530712}n^9}.
\end{align}
Then for all positive integers $n$
\begin{align}
\frac{1}{729^n}u(n)<\sum_{m=n}^{\infty}\frac{1}{m^{11}729^m}<
\frac{1}{729^n}v(n).
\end{align}
\end{lem}
\proof By manipulating \emph{Mathematica} software, it is not difficult to check
\begin{align}
&u(n)-\frac{1}{729}u(n+1)-\frac{1}{n^{11}}
=-\frac{\Psi_1(10;n)}{n^{11}(1+n)^{10}(11+728n)(739+728n)}<0,\\
&v(n)-\frac{1}{729}v(n+1)-\frac{1}{n^{11}}\\
=&\frac{\Psi_2(10;n)}{n^{11}(1+n)^9((529984 n+537992)(n-1)+ 489933)(489933+1067976n+529984 n^2)}\nonumber\\
>&0.\nonumber
\end{align}
 By multiplying the above inequalities by $\frac{1}{729^m}$, then adding these telescoping estimates from $m=n$ to $m=\infty$, we can finish the proof of Lemma 10.

\noindent{{\emph{The proof of Theorem 5.}} Just as the proof of Theorem 1, Theorem 5 can be proved similarly by Lemma 9 and 10. Here we omit the detail.\qed

\begin{lem} Let $g(m)=\mathrm{MC}_3(m)-\frac{1}{729}\mathrm{MC}_3(m+1)
-\varrho(m)$. We have for $m\ge 1$
\begin{align}
-\frac{D_4}{(m+\frac 12)^{12}}<g(m)+\frac{C_3}{m^{11}}
<-\frac{D_4}{m^{12}},
\end{align}
where
\begin{align*}
D_4=-&124280353667510106220979748750667909695624573786666800\\
&114069359390500727018449872352585153675322995125291007697/\\
&741003093304537143754634300103897398471054081282472499619\\
&4805971487253
427817067830381380462371074531080221491200.
\end{align*}
\end{lem}
\proof We can check by using \emph{Mathematica} software
\begin{align*}
&g(m)+\frac{C_3}{m^{11}}+\frac{D_4}{(m+\frac 12)^{12}}\\
=&\frac{\Psi_3(28;m)}
{2315\cdots 1600\cdot m^{11}(1+2m)^{12}(1+3m)^2(1+6m)^2(1+12m)^2(5+12m)^2\Psi_4(10;m)}
>0,
\end{align*}
\begin{align*}
&g(m)+\frac{C_3}{m^{11}}+\frac{D_4}{m^{12}}\\
=&-\frac{\Psi_5(17;m)}{7410\cdots 1200\cdot
m^{12}(1+3m)^2(1+6m)^2(1+12m)^2(5+12m)^2\Psi_6(10;m)}<0,
\end{align*}
and this completes the proof of the lemma.\qed
\begin{lem} Let
\begin{align}
&u_1(n)=\frac{1}{\frac{728}{729}n^{12}+\frac{4}{243}n^{11}}
,\quad v_1(n)=\frac{1}{\frac{728}{729}n^{12}},\\
&u_2(n)=\frac{1}{\frac{728}{729}(n+\frac 34)^{12}},
\quad v_2(n)=\frac{1}{\frac{728}{729}n^{12}+\frac{ 1460}{243}n^{11}}.
\end{align}
Then
\begin{align}
&\frac{1}{729^n}u_1(n)<\sum_{m=n}^{\infty}\frac{1}{m^{12}729^m}
<\frac{1}{729^n}v_1(n),\\
&\frac{1}{729^n}u_2(n)<\sum_{m=n}^{\infty}\frac{1}{(m+\frac 12)^{12}729^m}<\frac{1}{729^n}v_2(n).
\end{align}
\end{lem}
\proof The upper bound in the first inequalities is trivial.
By applying \emph{Mathematica} software, it isn't difficult to check
\begin{align}
&u_1(n)-\frac{1}{729}u_1(n+1)-\frac{1}{n^{12}}\\
=&-\frac{\Psi_1(11;n)}{4n^{12}(1+n)^{11}(3+182n)(185+182n)}<0,
\nonumber\\
&u_2(n)-\frac{1}{729}u_2(n+1)-\frac{1}{(n+\frac 12)^{12}}\\
=&-\frac{4096\Psi_2(23;n)}{91(1+2n)^{12}(3+4n)^{12}(7+4n)^{12}}<0,
\nonumber\\
&v_2(n)-\frac{1}{729}v_2(n+1)-\frac{1}{(n+\frac 12)^{12}}\\
=&\frac{\Psi_3(22;n)}
{4n^{11}(1+n)^{11}(1+2n)^{12}(1095+182n)(1277+182n)}>0.\nonumber
\end{align}
Now by multiplying the above inequalities by $\frac{1}{729^m}$, then adding these telescoping estimates from $m=n$ to $m=\infty$, we can finish the proof of Lemma 12.\qed

\begin{thm} For $n\ge 15$, we have the following inequalities
\begin{align}
\frac{27C_3}{364}\frac{1}{729^n(n+\frac 32)^{11}}<
E_3(n)<\frac{27C_3}{364}\frac{1}{729^n(n+\frac 12)^{11}}.
\end{align}
\end{thm}
\proof It follows from (5.3)
\begin{align}
E_3(n)=&\frac{2}{27}\left(\sum_{m=n}^{\infty}\frac{\varrho(m)}{729^m}
-\mathrm{MC}_3(n)\right).
\end{align}
By multiplying the inequalities (5.10)  by $\frac{1}{729^m}$, then adding these telescoping estimates from $m=n$ to $m=\infty$, we get
\begin{align}
&\frac{2}{27}\left(C_3\sum_{m=n}^{\infty}\frac{1}{729^mm^{11}}
+D_4\sum_{m=n}^{\infty}\frac{1}{729^mm^{12}}
\right)\\
<& E_3(n)<\frac{2}{27}\left(C_3\sum_{m=n}^{\infty}\frac{1}{729^mm^{11}}
+D_4\sum_{m=n}^{\infty}\frac{1}{729^m(m+\frac 12)^{12}}
\right).\nonumber
\end{align}
It follows from Lemma 10 and Lemma 12
\begin{align}
E_3(n)<\frac{2C_3}{27}\frac{1}{729^n}\left(
\frac{1}{\frac{728}{729}n^{11}+\frac{11}{729}n^{10}
-\frac{48059}{530712}n^9}
+\frac{D_4}{C_3}\frac{1}{\frac{728}{729}(n+\frac 34)^{12}}
\right).
\end{align}
By using \emph{Mathematica} software again, it is not difficult to verify for $n\ge 15$
\begin{align}
&\frac{1}{\frac{728}{729}n^{11}+\frac{11}{729}n^{10}
-\frac{48059}{530712}n^9}
+\frac{D_4}{C_3}\frac{1}{\frac{728}{729}(n+\frac 34)^{12}}
-\frac{1}{\frac{728}{729}(n+\frac 12)^{11}}\\
=&-\frac{1944\cdot \left(\Psi_4(21;n)(n-15)+2368\cdots 3535\right)}
{1418\cdots 2655\cdot n^9(1+2n)^{11}(3+4n)^{12}
\left((529984n+537992)(n-1)+489933\right)}<0.
\end{align}
This completes the proof of the right-hand side inequality in Theorem 6. Similarly, we have for $n\ge 12$
\begin{align}
E_3(n)>&\frac{2C_3}{27}\frac{1}{729^n}\left(
\frac{1}{\frac{728}{729}n^{11}+\frac{11}{729}n^{10}}
+\frac{D_4}{C_3}\frac{1}{\frac{728}{729}n^{12}}
\right)\\
>&\frac{27C_3}{364}\frac{1}{729^n(n+\frac 32)^{11}},\nonumber
\end{align}
here we use
\begin{align}
&\frac{1}{\frac{728}{729}n^{11}+\frac{11}{729}n^{10}}
+\frac{D_4}{C_3}\frac{1}{\frac{728}{729}n^{12}}
-\frac{1}{\frac{728}{729}(n+\frac 32)^{11}}\\
=&\frac{243\left(\Psi_5(11;n)(n-12)+2864\cdots 5343\right)}{3630\cdots 9680 n^{12}(3+2n)^{11}(11+728n)}>0.\nonumber
\end{align}
This finishes the proof of Theorem 6.
\section{Conclusions}
Our method may be used to establish similar results for many
series with a \emph{proper BBP-type term}. For example, such kind of series can be founded in~\cite{Bra,C-Z,Wei,Yee1,YK}. In this paper, we don't do any computation for these mathematical constants. However, for this question, we would like to pointed out that the computations of two main terms in our method (e.g. the second and third member of right hand side of (3.2) in Theorem 1 ) should play the same role, i.e. their computations should be ``matched".

In what follows, we give three examples to illustrate that the $k$-th correction function $MC_k(n)$ may be established occasionally
by a precise expression.

{\bf Example 1}  One has the following simple formula for Catalan constant
(see Entry 22 in~\cite{Wei})
\begin{align}
G=\frac 12\sum_{n=0}^{\infty}\frac{4^nn!^2}{(2n)!\cdot(2n+1)^2}.
\end{align}
By using \emph{Mathematica} software, one may check that the $k$-th correction $\mathrm{MC_k(n)}$ satisfies
\begin{align}
\mathrm{MC_0(n)}=\frac{\frac 12}{n+\frac 16},\quad
\mathrm{MC_k(n)}=\frac{\frac 12}{n+\frac 16+}\K_{j=1}^{k}\frac{a_j}{n+b_j}
, (k\ge 1),
\end{align}
where for $1\le k\le 20$
\begin{align}
a_k=\frac{2k^3(2k-1)^3}{(4k+1)(4k-1)^2(4k-3)},\quad
b_k=\frac{4k^2+2k-1}{2(4k-1)(4k+3)}.
\end{align}

{\bf Example 2} In 1668, Nicolas Mercator~\cite{Mer}  proved the following classical formula
\begin{align}
\ln 2=\sum_{n=1}^{\infty}\frac{1}{n 2^n}.
\end{align}
One may check that the $k$-th correction $\mathrm{MC_k(n)}$ has the form
\begin{align}
\mathrm{MC_0(n)}=\frac{2}{n+1},\quad
\mathrm{MC_k(n)}=\frac{2}{n+1+}\K_{j=1}^{k}\frac{a_j}{n+b_j}
, (k\ge 1),
\end{align}
where for $1\le k\le 20$
\begin{align}
a_k=-2\cdot k^2,\quad b_k=3k-2.
\end{align}

{\bf Example 3} For the series $\sum_{m=1}^{\infty}\frac{1}{(4m+1)^2}$, we can check
\begin{align}
\mathrm{MC_0(n)}=\frac{\frac{1}{16}}{n-\frac 14},\quad
\mathrm{MC_k(n)}=\frac{\frac{1}{16}}{n-\frac 14+}\K_{j=1}^{k}\frac{a_j}{n-\frac 14}
, (k\ge 1),
\end{align}
where for $1\le k\le 20$
\begin{align}
a_k=\frac{k^4}{4(2k-1)(2k+1)}.
\end{align}
Finally, we conjecture the above results should be true for all $k$.

\begin{flushleft}

Xiaodong Cao\\
Department of Mathematics and Physics, \\
Beijing Institute of Petro-Chemical Technology,\\
Beijing, 102617, P. R. China \\
 e-mail: caoxiaodong@bipt.edu.cn \\
\bigskip
Xu You \\
Department of Mathematics and Physics, \\
Beijing Institute of Petro-Chemical Technology,\\
Beijing, 102617, P. R. China \\
e-mail: youxu@bipt.edu.cn
\end{flushleft}

\end{document}